\documentclass[12pt]{article}
\usepackage{amsmath}
\usepackage{amscd}
\usepackage{mathrsfs}
\usepackage{amsfonts}
\usepackage{amsmath,amssymb}
\baselineskip 5.5pt \lineskip 5.5pt \numberwithin{equation}{section}

\def\ad{\mbox{ad}}

\def \Z{\hbox{$Z\hskip -5.2pt Z$}}

\def \C{\hbox{$C\hskip -5pt \vrule height 6pt depth 0pt \hskip 6pt$}}

\def\qed{\ \ \ifhmode\unskip\nobreak\fi\ifmmode\ifinner
         \else\hskip5pt\fi\fi
 \hbox{\hskip5pt\vrule width4pt height6pt depth1.5pt\hskip 1 pt}}
\def\a{\alpha}

\def\b{\beta}
\def\d{\delta}
\def\D{\Delta}
\def\g{\gamma}
\def\gi{\mathfrak{g}}

\def\L{\mathcal{L}}
\def\l{\lambda}

\def\cl{\centerline}

\def\Der{\mathrm{Der}}

\def\D{\Delta}

\def\vs{\vspace*}

\def\C{\mathbb{C}}

\def\Z{\mathbb{Z}}

\newtheorem{theo}{Theorem}[section]
\newtheorem{lemm}[theo]{Lemma}

\begin{document}
\cl {{\large\bf
 \vs{10pt} 2-local derivations on the twisted Heisenberg-Virasoro algebra}
\noindent\footnote{Supported by the National Science Foundation of
China (Nos. 11047030 and 11771122).
 }} \vs{6pt}

\cl{Yufang Zhao, Yongsheng Cheng}

\cl{ \small School
of Mathematics and Statistics, Henan
University, Kaifeng 475004, China} \vs{6pt}

\vs{6pt}

{\small
\parskip .005 truein
\baselineskip 10pt \lineskip 10pt

\noindent{{\bf Abstract:}\,2-local derivation is a generalized derivation for a Lie algebra,
which plays an important role to the study of local properties of the structure of the Lie algebra.
In this paper, we prove that every 2-local derivation on the twisted Heisenberg-Virasoro algebra
is a derivation.}
 \vs{5pt}

}

\noindent{\bf MR(2000) Subject Classification} 16E40, 17B56, 17B68.
\parskip .001 truein\baselineskip 8pt \lineskip 8pt

\vs{6pt}
\par
\cl{\bf\S1. \ Introduction}
\setcounter{section}{1}\setcounter{theo}{0}\setcounter{equation}{0}

Nowadays the theory of operator algebras plays an important role both in mathematics and physics.
This is motivated by the fact that in terms of operator algebras, their states, representations,
and derivations one can describe and investigate properties of model
systems in the quantum field theory and statistical physics.

Let $g$ be an algebra and $\D$ a map of $g$ into itself.
If $\D$ is linear and satisfies the identity $\D(xy)=\D(x)y+x\D(y)$
for all $x,y\in g$, then we call $\D$ a derivation of $g$.
Each element $a\in g$ defines a linear derivation $\D_a$ on $g$ given
by $\D_a(x)=ax-xa, x\in g$. Such derivations $\D_a$ are said to be inner.
In \cite{[S]}, as a generalization of derivation, P. \v{S}emrl
introduced the notion of 2-local derivations on algebras.
$\D$ (not necessarily linear) is called a $2$-local derivation of $g$, if for every pair
of elements $x,y\in g$, there exists a derivations
$\D_{x,y}: g\rightarrow g$ (depending on $x, y$) such that
$\D_{x, y}(x)=\D(x)$ and $\D_{x, y}(y)=\D(y)$.
For a 2-local derivation
on $\L$ and $k\in\C$, $x\in\L$, we have
$$\D(kx)=\D_{x,kx}(kx)=k\D_{x,kx}(x)=k\D(x).$$

Investigation of 2-local derivations on finite dimensional Lie algebras and infinite
dimensional Lie (super) algebras were initiated in papers \cite{[AK1], [AK2], [DGL], [T], [ZCZ]}.
In \cite{[AK1]}, the authors proved that every 2-local derivation on a
semi-simple Lie algebra is a derivation and that each finite-dimensional
nilpotent Lie algebra with dimension larger than two admits 2-local derivation
which is not a derivation. In \cite{[AK2], [DGL], [T], [ZCZ]},
the authors proved that 2-local derivations on the Witt algebra, super
Virasoro algebra, W-algebra $W(2, 2)$ and its superalgebra are derivations
and there are 2-local derivations on the so-called thin Lie algebra which are not derivations.

In this paper, we will study 2-local derivations on the twisted Heisenberg-Virasoro algebra.
The plan of this paper is as follows. In section 2, we give some preliminaries concerning
the twisted Heisenberg-Virasoro algebra. In section 3, we prove that
every 2-local derivation on the twisted Heisenberg-Virasoro algebra is automatically a derivation.
\vskip7pt

\cl{\bf\S2.\ Notations and Preliminaries}
\setcounter{section}{2}\setcounter{theo}{0}\setcounter{equation}{0}

Now let us recall the twisted Heisenberg-Virasoro algebra at level zero.
The twisted Heisenberg-Virasoro Lie algebra at level zero was first introduced in \cite{[ADKP]},
which is the universal central extension of the Lie algebra of differential operators on a circle of order at most one.
The twisted Heisenberg-Virasoro algebra $\L$ has a basis $$\{L_n,I_n,C_L,C_LI,C_I|n\in\Z\}$$
with the following commutation relations
$$\begin{array}{cc}
[L_n,L_m]=(n-m)L_{n+m}+\d_{n,-m}\frac{n^3-n}{12}C_L,\\[6pt]
[L_n,I_m]=-mI_{n+m}+\d_{n,-m}(n^2+n)C_{LI},\\[6pt]
[I_n,I_m]=n\d_{n,-m}C_I,\\[6pt]
[\L,C_L]=[\L,C_{LI}]=[\L,C_I]=0.
\end{array}$$
This Lie algebra has an infinite-dimensional Heisenberg subalgebra
and a Virasoro subalgebra. These subalgebras, however,
do not form a semidirect product, but instead, the natural action of the
Virasoro subalgebra on the Heisenberg subalgebra is twisted with a 2-cocycle.
The structure and representation theory of the twisted Heisenberg-Virasoro algebra were
studied by V. Kac, Y. Billing, C. Jiang, R. Shen and D.Liu etc (see \cite{[ADKP], [JJ], [LZ], [SJ]})

The following lemma comes from \cite{[SJ]}, which determines the derivations of the twisted Heisenberg-Virasoro
algebra.
\begin{lemm}\label{daozi}
$\Der\L=\ad\L\oplus\C D_1\oplus\C D_2\oplus\C D_2 $, where
\begin{align*}
D_1:&D_1|_{L\oplus \C C_L}=0,\ \ \ \ D_1|_{I\oplus\C C_{LI}}=id,\ \ \ \
     D_1C_1=2C_1;\\
D_2:&D_2(L_n)=nI_n+\d_{n,0}C_{LI},\ \ \ \ D_2I_n=-\d_{n,0}C_I,\ \ \ \
     D_2C_L=24C_{LI},\\
    &D_2C_{LI}=-C_I,\ \ \ \ D_2C_I=0;\\
D_3:&D_3(L_n)=(n+1)I_n,\ \ \ \ D_3C_L=24C_{LI},\ \ \ \ D_3C_{LI}=-C_I,
     \ \ \ \ D_3|_{I\oplus\C C_I}=0,
\end{align*}
where $n\in\Z,L=span_{\C}\{L_n,n\in\Z\},I=span_{\C}\{I_n,n\in\Z\}$
\end{lemm}

By Lemma \ref{daozi}, we can easily obtain
\begin{lemm}\label{xingshi}
Let $\D$ be a 2-local derivation on the twisted Heisenberg-Virasoro
algebra. Then for every $x,y\in\L$, there exists a derivation
$\D_{x,y}$ of $\L$ for which $\D_{x,y}(x)=\D(x)$ and
$\D_{x,y}(y)=\D(y)$ and it can be written as
\begin{align*}
\D_{x,y}&=\ad (\sum_{i\in\Z}(a_i(x,y)L_i
+b_i(x,y)I_i+l_1(x,y)C_L+l_2(x,y)C_{LI}+l_3(x,y)C_I)\\
&+\a(x,y)D_1+\b(x,y)D_2+\g(x,y)D_3.
\end{align*}
where $a_i, i\in \mathbb{Z}$, $b_i$, $\a$, $\b$, $\g$ are complex-valued
functions on $\L\times\L$ and $D_k$ for $k=1, 2, 3$ are given in Lemma \ref{daozi}.
\end{lemm}
\vskip7pt

\cl{\bf\S3.\ 2-local derivations on the twisted Heisenberg-Virasoro
algebra}
\setcounter{section}{3}\setcounter{theo}{0}\setcounter{equation}{0}
In this section, we will determine all 2-local derivations on the twisted Heisenberg-Virasoro
algebra.

\begin{lemm}\label{LiI0}
Let $\D$ be a 2-local derivation on $\gi$. For any but fixed $x\in\gi$.

(1) If $\D(L_i)=0$ for any $i\in \mathbb{Z}$, then
\begin{align*}
\D_{L_i,x}&=\ad(a_i(L_i,x)L_i+b_0(L_i,x)I_0
+l_1(L_i,x)C_L+l_2(L_i,x)C_{LI}+l_3(L_i,x)C_I)\\
&+\a(L_i,x)D_1+(1-\delta_{i, 0})\b(L_i,x)(D_2-\frac{i}{i+1}D_3), \ i\in \mathbb{Z}.
\end{align*}

(2) If $\D(I_0)=0$, then
\begin{align*}
\D_{I_0,x}&=\ad (\sum_{i\in\Z}(a_i(I_0,x)L_i
+b_i(I_0,x)I_i+l_1(I_0,x)C_L+l_2(I_0,x)C_{LI}+l_3(I_0,x)C_I)\\
&+\g(I_0,x)D_3.
\end{align*}
\end{lemm}
\noindent{\it Proof.~}
By Lemma \ref{xingshi}, for $x\in\{L_i,I_0\}$, we have
\begin{align*}
\D_{x,y}&=\ad (\sum_{i\in\Z}(a_i(x,y)L_i
+b_i(x,y)I_i+l_1(x,y)C_L+l_2(x,y)C_{LI}+l_3(x,y)C_I)\\
&+\a(x,y)D_1+\b(x,y)D_2+\g(x,y)D_3,
\end{align*}
where $a_i$, $b_i$, $\a$, $\b$, $\g$ are complex-valued
functions on $\L\times\L$ and $D_k, k=1, 2, 3$ are given in Lemma \ref{daozi}.

(1) If $\D(L_i)=0$, we have
\begin{align*}
\D(L_i)&=\D_{L_i,x}(L_i)\\
=&[\sum_{j\in\Z}(a_j(L_i,x)L_j+b_j(L_i,x)I_j+l_1(L_i,x)C_L
+l_2(L_i,x)C_{LI}+l_3(L_i,x)C_I,L_i]\\
&+\a(L_i,x)D_1(L_i)+\b(L_i,x)D_2(L_i)+\g(L_i,x)D_3(L_i)\\
=&\sum_{j\in\Z}((j-i)a_j(L_i,x)L_{i+j}+jb_j(L_i,x)I_{i+j})
-\frac{i^3-i}{12}a_{-i}(L_i,x)C_L\\
&+(i^2-i)b_{-i}(L_i,x)C_{LI}
+\b(L_i,x)(iI_i+\d_{i,0}C_{LI})+(i+1)\g(L_i,x)I_i\\
=&0.
\end{align*}
Thus we obtain $a_j(L_i,x)=0$ for $j\neq i$, $b_j(L_i,x)=0$ for $j\neq0$,
$\g(L_i,x)=\frac{i}{i+1}\b(L_i,x)$ for $j\neq0$ and
$\g(L_0,x)=\b(L_0,x)=0$.

(2) If $\D(I_0)=0$, we have
\begin{align*}
  \D(I_0)=&\D_{I_0,x}(I_0)\\
=&[\sum_{j\in\Z}(a_j(I_0,x)L_j+b_j(I_0,x)I_j+l_1(I_0,x)C_L
+l_2(I_0,x)C_{LI}+l_3(I_0,x)C_I,I_0]\\
&+\a(I_0,x)D_1(I_0)+\b(I_0,x)D_2(I_0)+\g(I_0,x)D_3(I_0)\\
=&\a(I_0,x)I_0-\b(I_0,x)C_I\\
=&0.
\end{align*}
Thus we obtain $\a(I_0,x)=\b(I_0,x)=0$.
\hfill$\Box$\vskip7pt

\begin{lemm}\label{L}
Let $\D$ be a 2-local derivation on $\gi$ such that $\D(L_0)=\D(L_1)=0$,
 then $\D(L_i)=0$, for any $i\in \mathbb{Z}$.
\end{lemm}
\noindent{\it Proof.~}
Since $\D(L_0)=\D(L_1)=0$, then we can assume that
\begin{align}
\D_{L_k, x}=&\ad(a_k(L_k, x)L_k+b_0(L_k, x)I_0
+l_1(L_k,x)C_L+l_2(L_k,x)C_{LI}\nonumber\\
 &+l_3(L_k,x)C_I)+\a(L_k,x)D_1
+(1-\delta_{k, 0})\b(L_1,x)(D_2-\frac{1}{2}D_3), \label{L1}
\end{align}
where $a_k$, $b_0$, $\a$, $\b$, $l_i, i=1, 2, 3$ are complex-valued
functions on $\L\times\L$, $k=0, 1$.
Then take $x=L_i$ in (\ref{L1}) for any $i$, we have
\begin{align*}
\D(L_i)=&\D_{L_0,L_i}(L_i)\\
=&[a_0(L_0,L_i)L_0+b_0(L_0,L_i)I_0+l_1(L_0,L_i)C_L+l_2(L_0,L_i)C_{LI}
\\
&+l_3(L_0,L_i)C_I,L_i]+\a(L_0,L_i)D_1(L_i)\\
=&-ia_0(L_0,L_i)L_i
\end{align*}
and
\begin{align*}
\D(L_i)&=\D_{L_1,L_i}(L_i)\\
&=[a_1(L_1,L_i)L_1+b_0(L_1,L_i)I_0+l_1(L_1,L_i)C_L+l_2(L_1,L_i)C_{LI}
\\
&+l_3(L_1,L_i)C_I,L_i]+\a(L_1,L_i)D_1(L_i)+\b(L_1,L_i)(D_2(L_i)-\frac{1}{2}D_3(L_i))\\
&=(1-i)a_1(L_1,L_i)L_{i+1}+\b(L_1,L_i)(iI_i-\frac{1}{2}(i+1)I_i)\\
&=(1-i)a_1(L_1,L_i)L_{i+1}-\frac{1}{2}(i-1)\b(L_1,L_i)I_i.
\end{align*}
Comparing the coefficients of the above equations, we get
$a_0(L_0,L_i)=a_1(L_1,L_i)=\b(L_1,L_i)=0$. It concludes that
$\D(L_i)=0$.
\hfill$\Box$\vskip7pt

\begin{lemm}\label{x}
Let $\D$ be a 2-local derivation on $\L$ such that $\D(L_i)=0$,
 then for any
 $x=\sum_{t\in\Z}(\a_tL_t+\b_tI_t)+k_1C_L+k_2C_{LI}+k_3C_I\in\L$, we
 have
 $$\D(x)=\l_x(\sum_{t\in\Z}\b_tL_t+k_2C_{LI}+k_3C_I), $$
 where $\l_x$ is a complex number depending on $x$.
\end{lemm}
\noindent{\it Proof.~}
For $x=\sum_{t\in\Z}(\a_tL_t+\b_tI_t)+k_1C_L+k_2C_{LI}+k_3C_I\in\L$,
since $\D(L_i)=0$, by Lemma \ref{L}, we have
\begin{align*}
\D(x)=&\D_{L_0,x}(x)\\
=&[a_0(L_0,x)L_0+b_0(L_0,x)I_0+l_1(L_0,x)C_L\\
&+l_2(L_0,x)C_{LI}
+l_3(L_0,x)C_I,x]+\a(L_0,x)D_1(x)\\
=&-\sum_{t\in\Z}ta_0(L_0,x)\a_tL_t
+\a(L_0,x)(\sum_{t\in\Z}\b_tI_t+k_2C_{LI}+2k_3C_I),
\end{align*}
and
\begin{align*}
\D(x)=&\D_{L_i,x}(x)\\
=&[a_0(L_i,x)L_i+b_0(L_i,x)I_0+l_1(L_i,x)C_L+l_2(L_i,x)C_{LI}
+l_3(L_i,x)C_I,x]\\
&+\a(L_1,x)D_1(x)+\b(L_1,x)(D_2(x)-\frac{i}{i+1}D_3(x))\\
=&a_0(L_i,x)(\sum_{t\in\Z}((t-i)\a_{t}L_{i+t}-t\b_tI_{i+t})
-\frac{n^3-n}{12}\a_{-i}C_L-(i^2+i)\b_{-i}C_{LI})\\
  &+\a(L_0,L_i)(\sum_{t\in\Z}\b_tI_t+k_2C_{LI}+2k_3C_I)
  +\b(L_1,x)(\sum_{t\in\Z}tI_t+\a_0C_{LI}-\b_0C_I\\
  &+24k_1C_{LI}-k_2C_I
  -\frac{i}{i+1}(\sum_{t\in\Z}(t+1)I_t+24k_1C_{LI}-k_2C_I))\\
=&a_0(L_i,x)(\sum_{t\in\Z}((t-i)\a_{t}L_{i+t}-t\b_tI_{i+t})
-\frac{n^3-n}{12}\a_{-i}C_L-(i^2+i)\b_{-i}C_{LI})\\
  &+\a(L_i,x)(\sum_{t\in\Z}\b_tI_t+k_2C_{LI}+2k_3C_I)
  +\frac{1}{i+1}\b(L_i,x)(\sum_{t\in\Z}(t-i)I_t\\
  &+24k_1C_{LI}-k_2C_I)+\b(L_i,x)(\a_0C_{LI}-\b_0C_I).
\end{align*}
If there exists an element $\a_t\neq0$, then by taking enough different
$i\in\Z$ in above equations, we obtain
$$a_0(L_i,x)=0, i\in \mathbb{Z}. $$
Using $\b_t\a(L_0,x)=\b_t\a(L_i,x)
+\frac{1}{i+1}(t-i)\a_t\b(L_i,x)$, we obtain
$$\b(L_i,x)=0, \ \a(L_0,x)=\a(L_i,x). $$
By the arbitrary of $x$, denote $\l_x=\a(L_i,x)$,
then
$$\D(x)=\l_x(\sum_{t\in\Z}\b_tL_t+k_2C_{LI}+k_3C_I). $$
\hfill$\Box$\vskip7pt

\begin{lemm}\label{lemm34}
Let $\D$ be a 2-local derivation on $\L$ such that $\D(L_0)=\D(L_1)=\D(I_0)=0$.
Then for any $p\in \mathbb{Z}^{*}$ and $y\in\L$, we have
\begin{eqnarray} \D_{L_{2p}+I_{p},y}=\!\!\!\!\!\!\!\!&\ad((a_{2p}(L_{2p}+I_{p},y)(L_{2p}+I_{p})
 +b_0(L_{2p}+I_{p},y)I_0+l_1(L_{2p}+I_{p}, y)C_L\nonumber \\
 &+l_2(L_{2p}+I_{p}, y)C_{LI}+l_3(L_{2p}+I_{p}, y)C_I)
 +\b(L_{2p}+I_{p},y)(D_2-\frac{2p}{2p+1}D_3). \label{lem34}
\end{eqnarray}
\end{lemm}
\noindent{\it Proof.~}
Using $\D(L_0)=\D(L_1)=0$ and Lemma \ref{L}, Lemma \ref{x},
for any $p\in \mathbb{Z}^{*}$ and $y\in\L$, there
exists $\l_{L_{2p}+I_{p}}\in \mathbb{C}$ satisfying
$$\D(L_{2p}+I_{p})=\l_{L_{2p}+I_{p}}I_p, \ p\in \Z^*. $$
Using $\D(I_0)=0$ and Lemma \ref{LiI0}, we have
\begin{align*}
\D(L_{2p}+I_{p})&=\D_{I_0,L_{2p}+I_{p}}(L_{2p}+I_{p})\\
&=[\sum_{i\in\Z}(a_i(I_0,L_{2p}+I_{p})L_i+b_i(I_0,L_{2p}+I_{p})I_i
+l_1(I_0,L_{2p}+I_{p})C_L\\
&+l_2(I_0,L_{2p}+I_{p})C_{LI}
+l_3(I_0,L_{2p}+I_p)C_I,L_{2p}+I_p]+\g(I_0,L_{2p}+I_p)D_3(L_{2p}+I_p)\\
&=\sum_{i\in\Z}(a_i(I_0,L_{2p}+I_{p})((i-2p)L_{i+2p}-pI_{i+p})
+ib_i(I_0,L_{2p}+I_{p})I_{2p+i})\\
&+(p^2+p)(b_{-p}(I_0,L_{2p}+I_{p})-a_{-p}(I_0,L_{2p}+I_{p}))C_{LI}
+\frac{8p^3-2p}{12}a_{-2p}(I_0,L_{2p}+I_{p})C_L\\
&+pb_{-p}(I_0,L_{2p}+I_{p})C_I+\g(I_0,L_{2p}+I_p)(2p+1)I_{2p}.
\end{align*}
Comparing the coefficients of the above equations, we obtain that
$$a_i(I_0,L_{2p}+I_{p})=b_i(I_0,L_{2p}+I_{p})=0, \ i\in \mathbb{Z}. $$
Then we have
$\l_{L_{2p}+I_{p}}=0$, and hence
$$\D(L_{2p}+I_{p})=0. $$
For every $y\in\L$, by Lemma \ref{xingshi}, we assume that
\begin{align*}
\D_{L_{2p}+I_p,y}&=\ad (\sum_{i\in\Z}(a_i(L_{2p}+I_p,y)L_i
+b_i(L_{2p}+I_p,y)I_i+l_1(L_{2p}+I_p,y)C_L\\
&+l_2(L_{2p}+I_p,y)C_{LI}+l_3(L_{2p}+I_p,y)C_I)+\a(L_{2p}+I_p,y)D_1\\
&+\b(L_{2p}+I_p,y)D_2
+\g(L_{2p}+I_p,y)D_3.
\end{align*}
So we have
\begin{align*}
\D(L_{2p}+I_p)&=\D_{L_{2p}+I_p,y}(L_{2p}+I_p)\\
=&[\sum_{i\in\Z}(a_i(L_{2p}+I_p,y)L_i
+b_i(L_{2p}+I_p,y)I_i+l_1(L_{2p}+I_p,y)C_L+l_2(L_{2p}+I_p,y)C_{LI}\\
&+l_3(L_{2p}+I_p,y)C_I,L_{2p}+I_p]+\a(L_{2p}+I_p,y)D_1(L_{2p}+I_p)\\
&+\b(L_{2p}+I_p,y)D_2(L_{2p}+I_p)
+\g(L_{2p}+I_p,y)D_3(L_{2p}+I_p).\\
=&\sum_{i\in\Z}(a_i(L_{2p}+I_p,y)((i-2p)L_{i+2p}-pI_{i+p})
+ib_i(L_{2p}+I_p,y)I_{i+2p})\\
&+\frac{8p^3-2p}{12}a_{-2p}(L_{2p}+I_p,y)C_L
+(p^2+p)(b_{-p}(L_{2p}+I_p,y)\\
&-a_{-p}(L_{2p}+I_p,y))C_{LI}+pb_{-p}(L_{2p}+I_p,y))C_I+\a(L_{2p}+I_p,y)I_p\\
&+2p\b(L_{2p}+I_p,y)I_{2p}+(2p+1)\g(L_{2p}+I_p,y)I_{2p}\\
=&0.
\end{align*}
From this, we obtain that
\begin{align*}
\a(L_{2p}+I_p,y)&=a_i(L_{2p}+I_p,y)=0, \ k\neq2p, \\
b_i(L_{2p}+I_p,y)&\neq0, \ k\neq0, p, \\
a_{2p}(L_{2p}+I_p,y)&=b_p(L_{2p}+I_p,y), \\
\g(L_{2p}+I_p,y)&=\frac{2p}{2p+1}\b(L_{2p}+I_p,y).
\end{align*}
Thus we obtain (\ref{lem34}).
\hfill$\Box$\vskip7pt

\begin{lemm}\label{zong}
Let $\D$ be a 2-local derivation on $\L$ such that
$\D(L_0)=\D(L_1)=\D(I_0)=0$,
 then $\D(x)=0$ for any $x\in\L$.
\end{lemm}
\noindent{\it Proof.~}
For $x=\sum_{t\in\Z}(\a_tL_t+\b_tI_t)+k_1C_L+k_2C_{LI}+k_3C_I\in\L$, where
$\a_t,\b_t (t\in \mathbb{Z}), l_1,l_2,l_3\in\C$.
Since $\D(L_0)=\D(L_1)=0$, by Lemma \ref{L}, we obtain $\D(L_i)=0$.
Furthermore, using Lemma \ref{x}, we have
$$\D(x)=\l_x(\sum_{t\in\Z}\b_tI_t+k_2C_{LI}+k_3C_I), $$
where $\l_x\in\C$.
On the other hand, for any $p\in \mathbb{Z}^*$, by Lemma \ref{lemm34}, we have (\ref{lem34}).
Thus
\begin{align*}
\D(x)=&\D_{L_{2p}+I_{p},x}(x)\\
=&[(a_{2p}(L_{2p}+I_{p},x)(L_{2p}+I_{p})+b_0(L_{2p}+I_{p},x)I_0
+l_1(L_0,x)C_L+l_2(L_0,x)C_{LI}\\
&+l_3(L_0,y)C_I,x]+\b(L_{2p}+I_p,y)(D_2(x)-\frac{2p}{2p+1}D_3(x))\\
=&\sum_{t\in\Z}a_{2p}(L_{2p}+I_{p},x)((t-2p)\a_tL_{2p+t}-t\b_tI_{2p+t}
+(p-t)\a_tI_{p+t}\\
&+\frac{8p^3-2p}{12}\a_{-2p}C_L+(4p^2+2p)\b_{-2p}C_{LI}+(p^2+p)\a_{-p}C_{LI}+p\b_{-p}C_{I})
\\
&+\b(L_{2p}+I_p,y)(\sum_{t\in\Z}t\a_tI_t+\a_0C_{LI}-\b_0C_I+24k_1C_{LI}-k_2C_I\\
&-\frac{2p}{2p+1}(\sum_{t\in\Z}(t+1)\a_tI_t+24k_1C_{LI}-k_2C_I)).
\end{align*}
If $\a_t=k_1=0$ for $k\in \mathbb{Z}$, i.e. $x=\sum_{t\in\Z}\b_tI_t+k_2C_{LI}+k_3C_I$, then we have
\begin{align*}
\D(x)=&\l_x(\sum_{t\in\Z}\b_tI_t+k_2C_{LI}+k_3C_I)\\
=&\sum_{t\in\Z}a_{2p}(L_{2p}+I_{p},x)(-t\b_tI_{2p+t}
+(4p^2+2p)\b_{-2p}C_{LI}\\
&+p\b_{-p}C_{I})-\b(L_{2p}+I_p,y)(\b_0C_I
+k_2C_I+\frac{2p}{2p+1}k_2C_I).
\end{align*}
This implies, for any $p\in \mathbb{Z}$ and the above any but fixed $x\in \L$,
\begin{align*}
\l_x\b_{2p}=&0, \\
\l_x\b_{2p+t}=&-t\b_ta_{2p}(L_{2p}+I_{p},x), \\
k_2\l_x=&(4p^2+2p)\b_{-2p}a_{2p}(L_{2p}+I_{p},x), \\
k_2\l_x=&p\b_{-p}a_{2p}(L_{2p}+I_{p},x)
-\b(L_{2p}+I_p,y)\frac{4p+1}{2p+1}k_2.
\end{align*}
Let $p$, $k_2$, $k_3$ run
 all integers, we conclude that
 $$\l_x=a_{2p}(L_{2p}+I_{p},x)=\b(L_{2p}+I_p,y)=0. $$
 Hence $\D(x)=0$.

If $\b_t=k_2=k_3=0$ for $t\in\Z$ i.e. $x=\sum_{t\in\Z}\a_tL_t+k_2C_{LI}+k_3C_I$
then by Lemma \ref{x}, we obtain $\D(x)=0$.

If both $\b_t,k_2,k_3$ and $\a_t,k_1$ for $t\in\Z$ are not
zero sequences. we assume $\a_0\neq0$,
then we have
\begin{align*}a_{2p}(L_{2p}+I_{p},x)&=0, \\
\l_x\a_t=&0, t=2,3,\dots,2l, \\
\l_xk_4=&(-1)^{l+\frac{1}{2}}(2l)!a_0(e+p_0,x)\\
\l_x\b_0&=-\frac{2p}{2p+1}\b(L_{2p}+I_p,y).
\end{align*}
Let $p$ run all integers, we conclude that
$\l_x=\b(L_{2p}+I_p,y)=0$, then $\D(x)=0$.
\hfill$\Box$\vskip7pt

\begin{theo}
Every 2-local derivation on $\L$ is a derivation.
\end{theo}
\noindent{\it Proof.~}
Let $\D$ is a 2-local derivation on $\L$. There exists a derivation $\D_{L_0,L_1}$
such that
$$\D(L_0)=\D_{L_0,L_1}(L_0) \ \mbox{and} \ \D(L_1)=\D_{L_0,L_1}(L_1). $$
Denote $\D(1)=\D-\D_{L_0,L_1}$. Then $\D(1)$ is a 2-local derivations satisfying
$$\D(1)(L_0)=\D(1)(L_1)=0. $$
By Lemma \ref{L}, we have $\D(1)(L_i)=0$.
From this with Lemma \ref{x}, we have
$$\D(1)(I_0)=\l_{I_0}I_0, \ \mbox{where} \ \l_{I_0}\in\C. $$
Set $\D(2)=\D(1)-\l_{I_0}D_1$. Then $\D(2)$ is a 2-local
derivation such that
\begin{align*}
&\D(2)(L_0)=\D(1)(L_0)-\l_{I_0}D_1(L_0)=0,\\
&\D(2)(L_1)=\D(1)(L_1)-\l_{I_0}D_1(L_1)=0,\\
&\D(2)(I_0)=\D(1)(I_0)-\l_{I_0}D_1(I_0)=\l_{I_0}D_1(I_0)-\l_{I_0}D_1(I_0)=0.
\end{align*}
By lemma \ref{zong}, we have $\D(2)=\D-\D_{L_0,L_1}-\l_{I_0}D_1\equiv0$. Hence
$\D=\D_{L_0,L_1}+\l_{I_0}D_1$ is a derivation.
\hfill$\Box$

 \end{document}